\documentclass[english,aps,preprint]{revtex4}
\usepackage[latin9]{inputenc}
\usepackage{mathrsfs}
\usepackage{mathtools}
\usepackage{amsmath}
\usepackage{amssymb}
\PassOptionsToPackage{normalem}{ulem}
\usepackage{ulem}

\makeatletter
\usepackage[english]{babel}
\usepackage{graphicx}
\usepackage{bm}
\usepackage{dcolumn}
\usepackage{xcolor}
\usepackage{bbm}
\usepackage{soul}
\usepackage{verbatim}

\usepackage{babel}

\begin{document}

\global\long\def\thefootnote{\arabic{footnote}}
\setcounter{footnote}{0}

\title{\sc\Large{\vspace*{-0.5cm} Definite integral of a Laguerre polynomial and exponentials}
\vspace*{0.5cm}}

\author{D. Gomez Dumm$^1$ and N.~N.\ Scoccola$^{2,3}$ \vspace*{0.5cm}}

\affiliation{$^{1}$ IFLP, CONICET $-$ Departamento de F\'{\i}sica, Facultad de Ciencias Exactas,
Universidad Nacional de La Plata, C.C. 67, 1900 La Plata, Argentina}
\affiliation{\small $^{2}$ Departamento de F\'isica, Comisi\'{o}n Nacional de Energ\'{\i}a At\'{o}mica, \\
Avenida del Libertador 8250, 1429 Buenos Aires, Argentina}
\affiliation{$^{3}$ CONICET, Rivadavia 1917, 1033 Buenos Aires, Argentina
\vspace{1cm}}

\begin{abstract}
In our investigations on the effect of strong magnetic fields on the
properties of elementary particles we have been faced with a definite
integral of the form
$$\int_0^{2\pi}d\theta\ L_{n}(s^2+t^2+2st\cos\theta)\
e^{-ik\theta}\, \exp{(-st\,e^{i\theta})}\ , $$ where $L_n(x)$ is a Laguerre
polynomial, $s$ and $t$ are real numbers and $n$ and $k$ are integers, with
$n \geq 0$. In the present article we show that this integral can be solved
analytically. The result can be used to get an alternative proof of an
addition formula for Laguerre polynomials.
\end{abstract}
%\date{\today}

\maketitle
\newpage

%\renewcommand{\thefootnote}{\arabic{footnote}}
%\setcounter{footnote}{0}

%\section{Statement}
%\label{intro}

The study of the effect of a magnetic field on elementary particles, in the
framework of quantum field theory, leads to various integrals and sums that
involve Laguerre polynomials (see e.g.\
Refs.~\cite{Andersen:2014xxa,Miransky:2015ava}). In particular, we have been
faced~\cite{GomezDumm:2017jij} with the definite integral
\begin{equation}
I_{n,k}(s,t) \ \equiv \ \frac{1}{2\pi}\int_0^{2\pi}d\theta\
L_{n}(s^2+t^2+2st\cos\theta)\ e^{-ik\theta}\, \exp (-st\,e^{i\theta}) \ .
\label{integral}
\end{equation}
Here, $s$ and $t$ are real numbers, $n$ and $k$ are integers, and $L_n(x)$
is a Laguerre polynomial, with $n\geq 0$. Interestingly, it can be shown
analytically that the result of this integral is
\begin{equation}
I_{n,k}(s,t) \ = \ \left\{
\begin{array}{ll}
\displaystyle (-1)^k \frac{n!}{(n+k)!}\, (st)^k\, L_n^k(s^2)\,L_n^k(t^2) & \ \ {\rm if}\ \ k \geq -n\ \ , \\
 \ \ 0 & \ \ {\rm if}\ \ k < -n \ \ , \end{array}
\right.
\label{result}
\end{equation}
where $L_n^k(x)$ is a generalized Laguerre polynomial. Notice that the
result in Eq.~(\ref{result}) can be factorized into functions that depend
only on $s$ and $t$, a fact that might be very useful when dealing with
expressions where $I_{n,k}(s,t)$ appears. We conclude this work by showing
that our result can be used to get an alternative proof of an addition
formula for Laguerre polynomials.

\section*{Proof}

Let us consider separately the cases  $k<-n$ and $k\geq -n\,$. The result in
Eq.~(\ref{result}) for $k<-n$ can be proved by noticing that $I_{n,k}(s,t)$
can be written as an integral on the complex plane $z$ along a circle of
radius one. One has
\begin{equation}
I_{n,k}(s,t) \ = \ -\,\frac{i}{2\pi}\int_C dz\ z^{-(k+1)}\,
L_{n}\bigg(s^2+t^2+st\Big(z + \frac{1}{z}\Big)\bigg)\ e^{-stz}\ ,
\label{intz}
\end{equation}
where the integration contour $C$ goes along the circle in anticlockwise
direction. Since $L_n(x)$ is a polynomial of degree $n$, it is easy to see
that for $k+n\geq 0$ the integrand in the above expression is analytical in
all points of the complex plane except the origin, $z=0$, where it has pole
of order $k+n+1$. On the other hand, if $k< -n$ there are no poles; thus,
from Cauchy's theorem one immediately finds $I_{n,k}(s,t)=0$, as stated in
Eq.~(\ref{result}).

For $k\geq -n$, the result in Eq.~(\ref{result}) can be proved using
mathematical induction. To carry our this proof, let us previously
show the validity of the differential equation
\begin{equation}
s\,\frac{d I_{n+1,k}(s,t)}{ds} - \big[2(n+1)+k]\, I_{n+1,k}(s,t) \
= \ \big[2t^2- 2(n+1)-k\big]\, I_{n,k}(s,t) - t\,\frac{d
I_{n,k}(s,t)}{dt}\ ,
\label{diffeq1}
\end{equation}
where $n\geq 0$ and $s,t$ are assumed ---for the moment--- to be positive
real numbers. It is seen that Eq.~(\ref{diffeq1}) follows from a series of
relations that involve integrals of Laguerre polynomials. The following
properties of these polynomials, valid for all $k$ and $m$ (see Sec.~8.971
of Ref.~\cite{grad}), are to be taken into account:
\begin{eqnarray}
L_{m-1}^{k+1}(x) + L_m^k(x) & = & L_m^{k+1}(x) \ ,
\label{uno} \\ [0.1cm]
x L_m^{k+1}(x) & = & (m+k+1)\,L_m^k(x) - (m+1)L_{m+1}^k(x) \ ,
\label{dos} \\ [0.2cm]
\frac{d L_{m+1}^k(x)}{dx} & = & -L_m^{k+1}(x) \ .
\label{derivada}
\end{eqnarray}

We start by considering the expression of $I_{n,k}(s,t)$ given in
Eq.~(\ref{intz}), assuming $k\geq -n$. According to Cauchy's theorem, since
the integrand has a single pole in $z=0$, the result of the integral cannot
change if the circle size is rescaled by an arbitrary positive factor. In
particular, let us take a circle of radius $s$, and parameterize the path by
$z=s\exp(i\theta)$. We get
\begin{equation}
I_{n,k}(s,t) \ = \ \frac{1}{2\pi}\int_0^{2\pi} d\theta\ s^{-k}\,
L_n\bigg(s^2+t^2+st\Big(s\,e^{i\theta} + \frac{1}{s\,e^{i\theta}}\Big)\bigg)\,
e^{-ik\theta}\, \exp{(-s^2t\,e^{i\theta})}\ .
\end{equation}
Next we can take the derivative with respect to $s$. We find
\begin{eqnarray}
\frac{d I_{n,k}(s,t)}{ds} & = & -\,\frac{1}{2\pi}\int_0^{2\pi} d\theta\ s^{-k}\,
e^{-ik\theta}\,\exp{(-s^2t\,e^{i\theta})} \nonumber \\
& & \times\,\bigg[\Big(\frac{k}{s}+2st\,e^{i\theta}\Big)\,
L_n\bigg(s^2+t^2+st\Big(s\,e^{i\theta} + \frac{1}{s\,e^{i\theta}}\Big)\bigg)\nonumber \\
& & \ \ \ \ +\,2s\Big(1 + t\,e^{i\theta}\Big)\,L_{n-1}^1
\bigg(s^2+t^2+st\Big(s\,e^{i\theta} + \frac{1}{s\,e^{i\theta}}\Big)\bigg)\bigg]\ ,
\end{eqnarray}
where we have made use of the relation in Eq.~(\ref{derivada}).
Once again, this last expression can be written as an integral on
the complex plane along a circle of radius $s$, with
$z=s\exp(i\theta)$. The integrand is shown to have, as in
Eq.~(\ref{intz}), a single pole of order $k+n+1$ in $z=0$. Thus,
we can rescale back the circle size to radius one. In this
way we get a first relation,
\begin{eqnarray}
s\,\frac{d I_{n,k}(s,t)}{ds} & = & -k I_{n,k}(s,t)-2s^2
I^1_{n-1,k}(s,t)- 2 st
\Big[I_{n,k-1}(s,t)+I^1_{n-1,k-1}(s,t)\Big]\ ,
\label{derivada1}
\end{eqnarray}
where we have introduced the function $I_{n,k}^1(s,t)$. The latter is
defined as an integral like that in Eq.~(\ref{integral}), in which we
replace $L_n(s^2+t^2+2st\cos\theta) \to L_n^1(s^2+t^2+2st\cos\theta)$.

We can follow a similar procedure (rescaling, derivative with
respect to $s$, inverse rescaling) considering a circle of radius
$1/s$. This leads to our second relation,
\begin{equation}
s\,\frac{d I_{n,k}(s,t)}{ds} \ = \ k I_{n,k}(s,t)-2s^2
I^1_{n-1,k}(s,t)- 2 st I^1_{n-1,k+1}(s,t)\ .
\label{derivada2}
\end{equation}
In addition, one can take a derivative with respect to $s$ directly on
Eq.~(\ref{integral}). One finds in this way
\begin{equation}
s\,\frac{d I_{n,k}(s,t)}{ds} \ = \ -2s^2 I^1_{n-1,k}(s,t)-st
\Big[I_{n,k-1}(s,t)+I^1_{n-1,k-1}(s,t)+ I^1_{n-1,k+1}(s,t)\Big]\ .
\end{equation}
Noticing that both the functions $I_{n,k}(s,t)$ and $I^1_{n,k}(s,t)$ are
symmetric under the exchange $s\leftrightarrow t$, from this last equation
one gets a third useful relation,
\begin{equation}
s\,\frac{d I_{n,k}(s,t)}{ds} - t\,\frac{d I_{n,k}(s,t)}{dt}\ = \ -2(s^2
- t^2) I^1_{n-1,k}(s,t)\ .
\label{derivada3}
\end{equation}

Our goal is to obtain a differential equation in which we can get
rid of the integrals over the extended Laguerre polynomials
$L^1_n(x)$. This can be done with the help of two more relations
that can be obtained using the properties of Laguerre polynomials.
From Eq.~(\ref{uno}), it is seen that the integrals satisfy a
fourth relation,
\begin{equation}
I^1_{n-1,k}(s,t) + I_{n,k}(s,t) \ = \ I^1_{n,k}(s,t)\ .
\label{int_uno}
\end{equation}
In addition, let us consider the sum
\begin{eqnarray}
(s^2+t^2)\, I^1_{n,k}(s,t) + st\Big[I^1_{n,k-1}(s,t) + I^1_{n,k+1}(s,t)\Big]\ =
& & \nonumber \\
& & \hspace{-8cm} \frac{1}{2\pi}\int_0^{2\pi}d\theta\
(s^2+t^2+2st\cos\theta)\,
L_{n}^1(s^2+t^2+2st\cos\theta)\ e^{-ik\theta}\, \exp{(-st\,e^{i\theta})} \ .
\end{eqnarray}
Using Eq.~(\ref{dos}), it is easy to show that the right hand side
can be expressed in terms of integrals of the form
(\ref{integral}). One gets in this way a fifth relation,
\begin{equation}
(s^2+t^2)\, I^1_{n,k}(s,t) + st\Big[I^1_{n,k-1}(s,t) + I^1_{n,k+1}(s,t)\Big]
\ = \ (n+1) \Big[I_{n,k}(s,t) - I_{n+1,k}(s,t)\Big]\ .
\label{int_dos}
\end{equation}

Now, taking into account the relations in
Eqs.~(\ref{derivada1}), (\ref{derivada3}), (\ref{int_uno}) and (\ref{int_dos}),
one obtains
\begin{eqnarray}
2 s^2 I^1_{n,k}(s,t) + 2 st I^1_{n,k+1}(s,t) & = &
\big[2(n+1)+k-2t^2\big]\, I_{n,k}(s,t) \nonumber \\
& & - 2(n+1) I_{n+1,k}(s,t) +
t\,\frac{d I_{n,k}(s,t)}{dt}\ .
\label{last}
\end{eqnarray}
Finally, considering Eq.~(\ref{derivada2}) with the replacement $n\to
n+1$, it is seen that the left hand side of Eq.~(\ref{last}) can be written
in terms of $I_{n+1,k}(s,t)$. In this way one arrives at the differential
equation (\ref{diffeq1}).

\hfill

We proceed now with the proof of our result in Eq.~(\ref{result}) for $k\geq
-n$. As mentioned above, this will be done using mathematical induction on
$n$.

Let $n=0$. After a series expansion of the last exponential in the integrand
in Eq.~(\ref{integral}) one gets
\begin{equation}
I_{0,k}(s,t) \ = \ \frac{1}{2\pi}\sum_{\ell=0}^\infty(-1)^\ell\,
\frac{(st)^\ell}{\ell!}\int_0^{2\pi}d\theta\
e^{i(\ell-k)\theta}\ .
\end{equation}
The integral on the right hand side is equal to $2\pi\delta_{\ell k}\,$. Since
$L_0^k(x)=1$, for $k\geq 0$ this leads to
\begin{equation}
I_{0,k}(s,t) \ = \ (-1)^k\, \frac{(st)^k}{k!}\ ,
\end{equation}
consistently with Eq.~(\ref{result}).

Next, let us assume that the result in Eq.~(\ref{result}) is valid
for a given $n$. We have to show that then it is also valid for
$n+1$. For this purpose we take into account the relation
in Eq.~(\ref{diffeq1}), replacing $I_{n,k}(s,t)$ in the right hand
side by the expression in Eq.~(\ref{result}). Using once again
the properties in Eqs.~(\ref{uno}-\ref{derivada}), we get
\begin{equation}
s\,\frac{d I_{n+1,k}(s,t)}{ds} - \big[2(n+1)+k]\, I_{n+1,k}(s,t) \ = \
h(t)\, s^k\, L_n^k(s^2)\ ,
\label{diffeq2}
\end{equation}
where the function $h(t)$ is given by
\begin{equation}
h(t) \ = \ 2(-1)^{k+1} \frac{(n+1)!}{(n+k)!}\; t^k\, L_{n+1}^k(t^2)\ .
\end{equation}
We have to prove that the expression for $I_{n+1,k}(s,t)$ arising from
Eq.~(\ref{result}) with $n\to n+1$ is the most general solution of
Eq.~(\ref{diffeq2}). Thus, we propose to write $I_{n+1,k}(s,t)$ in the
(general) form
\begin{equation}
I_{n+1,k}(s,t) \ = \ -\,\frac{h(t)}{2(n+k+1)}\; s^k\,L^k_{n+1}(s^2) + F(s,t)\ .
\end{equation}
Replacing in Eq.~(\ref{diffeq2}) and using the properties in
Eqs.~(\ref{dos}) and (\ref{derivada}) it is easy to show that most terms
cancel, in such a way that $F(s,t)$ turns out to be a solution of the
homogeneous equation
\begin{equation}
s\,\frac{d F(s,t)}{ds} - \big[2(n+1)+k]\, F(s,t) \ = \ 0\ .
\end{equation}
Taking into account the symmetry $s\leftrightarrow t$, it is seen that the
solution of this last equation has to be of the form $F(s,t) = a
(st)^{2(n+1)+k}$, where $a$ is a constant. Thus, we have
\begin{equation}
I_{n+1,k}(s,t) \ = \ (-1)^k \frac{(n+1)!}{(n+k+1)!}\, (st)^k\,
L_{n+1}^k(s^2)\,L_{n+1}^k(t^2) + a\,(st)^{2(n+1)+k}\ .
\label{result_a}
\end{equation}
This is the desired result, except for the presence of an unwanted
additional term on the right hand side. To prove that this term vanishes
(i.e., that $a=0$), let us once again use the strategy of writing
$I_{n+1,k}(s,t)$ as an integral along a circle in the complex plane, just
as in Eq.~(\ref{intz}), with $n\to n+1$. If we now rescale the circle radius
from $r=1$ to $r=1/(st)$, and use the parametrization $z=e^{i\theta}/(st)$,
we obtain
\begin{equation}
I_{n+1,k}(s,t) \ = \ \frac{1}{2\pi}\int_0^{2\pi} d\theta\ (st)^k\,
L_{n+1}\Big(s^2+t^2+e^{i\theta} + (st)^2\,e^{-i\theta}\Big)\,
e^{-ik\theta}\, \exp(-e^{i\theta})\ .
\end{equation}
It is explicitly seen that, as expected, the function $I_{n+1,k}(s,t)$ is a
polynomial in $s^2$, $t^2$ and $st\,$. In particular, the coefficient of the
term of order $(st)^{2(n+1)+k}$ is given by
\begin{eqnarray}
C_{2(n+1)+k} & = & \frac{(-1)^{n+1}}{2\pi(n+1)!}\,\int_0^{2\pi}d\theta\
e^{-i(k+n+1)\theta}\,\exp(-e^{i\theta}) \nonumber \\
& = & \frac{(-1)^k}{(n+1)!(k+n+1)!}\ ,
\label{coef}
\end{eqnarray}
where we have taken into account the fact that for a Laguerre polynomial
$L_m(x)$ the coefficient of the term of degree $m$ is $(-1)^m/m!\,$. The
result in Eq.~(\ref{coef}) can be compared with the corresponding
coefficient that arises from the right hand side of Eq.~(\ref{result_a}),
namely
\begin{eqnarray}
C_{2(n+1)+k} & = & (-1)^{k}\,\frac{(n+1)!}{(n+k+1)!}\,
\left[\frac{(-1)^{n+1}}{(n+1)!}\right]^2 \, +\, a
 \nonumber \\
& = & \frac{(-1)^k}{(n+1)!(k+n+1)!}\,+\,a \ .
\label{coef2}
\end{eqnarray}
The agreement of Eqs.~(\ref{coef}) and (\ref{coef2}) implies $a=0$. In this
way, we have shown that if one assumes that the result in
Eq.~(\ref{result}) is valid for a given $n\geq 0$, then it has to be
also valid for $n+1$.

We recall that we have made the assumption of taking $s$ and $t$ to be
positive real numbers. Now, it can be easily seen that the result in
Eq.~(\ref{result}) can be extended to negative values of $s$ or $t$ by
changing $\theta\to\theta+\pi$ in the integral in Eq.~(\ref{integral}). The
validity of our result for $s=0$ and/or $t=0$ can also be trivially
verified. This completes our proof.

\hfill

It is worth noting that, as usual, the analytical result obtained here can
be used to get other mathematical expressions. As an example, the integrand
in Eq.~(\ref{integral}) can be written in terms of a series of Bessel
functions using the relation~\cite{watson}
\begin{equation}
e^{- is t\sin\theta} \ = \ J_0(st) \ +\, \sum_{\ell\,=1}^\infty
\ (-1)^\ell \,\Big[ e^{i\ell\,\theta} + (-1)^\ell\,e^{-i\ell\,\theta} \Big]
\, J_\ell(st)\ .
\label{besselcos}
\end{equation}
Then, one can multiply the integral by a Bessel function $J_k(st)$ and sum
over $k$. With the aid of the relation~\cite{watson}
\begin{equation}
\sum_{k=-\infty}^\infty (-1)^k\,J_{\ell-k}(x)\,J_k(x) \ = \ \delta_{\ell\, 0}\ ,
\end{equation}
it is easy to arrive at the expression
\begin{eqnarray}
\frac{1}{2\pi}\int_0^{2\pi}d\theta\;
e^{-st\cos\theta}\,L_n(s^2+t^2+2 s\, t\cos\theta) & = & \nonumber \\
& & \hspace{-3cm}\frac{n!}{(s\, t)^n}\,\sum_{k=0}^\infty \frac{(s\,t)^k}{k!}\,
L_n^{k-n}(s^2)\,L_n^{k-n}(t^2)\,J_{k-n}(s\,t)\ ,
\label{byproduct}
\end{eqnarray}
which holds for nonzero real values of $s$ and $t$, and integer $n\geq 0$.
We believe that Eq.~(\ref{byproduct}) can be seen as an interesting
byproduct of our main result in Eq.~(\ref{result}).

\hfill

After the completion of our demonstration, we became aware of the existence
of an addition formula for Laguerre polynomials stated years ago by
H.\ Bateman~\cite{bateman} and generalized by T.\
Koornwinder~\cite{koornwinder}. This formula reads
\begin{equation}
\exp (xy\, e^{i\theta})\,L_{n}(x^2+y^2-2xy\cos\theta)\ = \
\sum_{j=0}^\infty\,
\big(xy\,e^{i\theta}\big)^{j-n}\,\frac{n!}{j!}\, L_n^{j-n}(x^2)\,L_n^{j-n}(y^2)
\ ,
\label{bateman}
\end{equation}
where $x$ and $y$ are real numbers. Using this expression one can obtain a
much simpler proof of the result in Eq.~(\ref{result}).
Setting $s=-x$, $t=y$, multiplying by $e^{-ik\theta}$, and integrating over
$\theta$ one gets~\cite{vanassche}
\begin{eqnarray}
\frac{1}{2\pi}\int_0^{2\pi}d\theta\
L_{n}(s^2+t^2+2st\cos\theta)\ e^{-ik\theta}\, \exp (-st\,e^{i\theta}) & = &
\nonumber \\
& & \hspace{-5.9cm}\frac{1}{2\pi}\sum_{j=0}^\infty\,
(-st)^{j-n}\,\frac{n!}{j!}\, L_n^{j-n}(x^2)\,L_n^{j-n}(y^2)\,
\int_0^{2\pi}d\theta\ e^{i\theta(-k+j-n)}\ .
\end{eqnarray}
The integral on the right hand side is zero for $j\neq k+n$ and $2\pi$ for
$j=k+n$. Therefore, for $k<-n$ all terms in the sum are zero, while for
$k\geq -n$ one gets our result in Eq.~(\ref{result}).

Although this last procedure is more immediate, we believe that the proof
presented in this work offers a different method that is worth to be taken
into account. Moreover, the problem can be put the other way round: our
demostration of Eq.~(\ref{result}) can be used to get an alternative proof
of Bateman's addition formula. As stated in Ref.~\cite{koornwinder}, proofs
of this formula have been previously given by Buchholz~\cite{buchholz},
Carlitz~\cite{carlitz} and Miller~\cite{miller}, following different
methods. In fact, starting from the expression in Eq.~(\ref{result}), we can
multiply by $e^{ik\phi}$ and sum over $k$ to obtain
\begin{eqnarray}
\sum_{k=-\infty}^\infty\,\frac{1}{2\pi}\int_0^{2\pi}d\theta\
L_{n}(s^2+t^2+2st\cos\theta)\ e^{ik(\phi-\theta)}\, \exp (-st\,e^{i\theta}) & = &
\nonumber \\
& & \hspace{-5.9cm}\sum_{k=-n}^\infty\, (-1)^k \frac{n!}{(n+k)!}\, (st)^k\,
L_n^k(s^2)\,L_n^k(t^2)\,e^{ik\phi}\ .
\label{delta}
\end{eqnarray}
Here the sum on the right hand side is started from $k=-n$, since ---as
proved above--- the integral is zero for lower values of $k$. Now, taking
into account the representation of the delta function
\begin{equation}
\frac{1}{2\pi}\,\sum_{k=-\infty}^\infty\ e^{ik(\phi-\theta)} \ = \
\delta(\phi-\theta)\ ,
\end{equation}
one can integrate over $\theta$ on the left hand side of Eq.~(\ref{delta})
to get
\begin{equation}
\exp (-st\, e^{i\phi})\,L_{n}(s^2+t^2+2st\cos\phi)\ = \
\sum_{k=-n}^\infty\,
\big(-st\,e^{i\phi}\big)^{k}\,\frac{n!}{(n+k)!}\, L_n^k(s^2)\,L_n^k(t^2)
\ ,
\end{equation}
which after variable changes gives indeed Bateman's addition formula,
Eq.~(\ref{bateman}). This proof represents an alternative to those proposed
in Refs.~\cite{buchholz,carlitz,miller}.

\begin{acknowledgments}
We are grateful to W.\ Van Assche, E.\ Koelink and T.\ Koornwinder for
useful remarks that made us become aware of the addition formula considered
in the final part of this article. This work has been partially funded by
CONICET (Argentina) under Grant No. PIP2022-GI-11220210100150CO, by ANPCyT
(Argentina) under Grant No.~PICT20-01847 and by the National University of
La Plata (Argentina), Project No.~X960.
\end{acknowledgments}

\end{document}